\documentclass[11pt, reqno]{amsart}
\usepackage{geometry}
\usepackage{layout}
\usepackage{amssymb}
\usepackage{amsmath}
\usepackage{amsthm}
\usepackage{braket}
\usepackage{pdfpages}
\usepackage{graphicx}
\usepackage{caption}
\usepackage{subcaption}
\usepackage{hyperref}
\usepackage{enumitem}
\usepackage{todonotes}

\numberwithin{equation}{section}
\theoremstyle{plain}
\newtheorem{theorem}{Theorem}[section]

\newtheorem{prop}[theorem]{Proposition}
\newtheorem{lem}[theorem]{Lemma}

\newtheorem{prob}[theorem]{Problem}

\theoremstyle{definition}
\newtheorem{defi}[theorem]{Definition}
\newtheorem{rmk}[theorem]{Remark}

\newcommand{\R}{\mathbb{R}}
\newcommand{\PP}{\mathbb{P}}
\newcommand{\RR}{\mathbb{R}}
\newcommand{\SSS}{\mathbb{S}}

\newcommand{\bv}{{\mathbf {v}}}

\newcommand{\sD}{{\mathcal D}}
\newcommand{\sE}{{\mathcal E}}
\newcommand{\sO}{{\mathcal O}}

\newcommand{\sDD}{{\mathcal  D_{\sig}}}

\newcommand{\shf}{{\mathfrak{s}}}

\newcommand{\sha}{{\mathfrak{s}}}

\newcommand{\sih}{{\mathfrak{s}}}
\newcommand{\sig}{{\sigma}}
\newcommand{\sh}{{\mathfrak{s}}}
\newcommand{\si}{\sh^\ast}
\newcommand{\shu}{\sh_u}
\newcommand{\shb}{\sh_b}
\newcommand{\siu}{\si_u}
\newcommand{\sib}{\si_b}
\newcommand{\ls}{\mathfrak{l}}
\newcommand{\cs}{\mathfrak{c}}

\newcommand{\pol}[1]{{#1}^\circ}

\newcommand{\sprod}[2]{\langle {#1} , {#2} \rangle} %

\newcommand{\defn}[1]{{\color{blue} \em #1}}

\numberwithin{figure}{section}

\title{Moser's Shadow Problem}
\author{Jeffrey C. Lagarias}
\address{Dept. of Mathematics, Univ. of Michigan, 530 Church Street, Ann Arbor, MI 48109-1043 USA}
\email{lagarias@umich.edu}
\author{Yusheng Luo}
\address{Dept. of Mathematics, Harvard University, One Oxford Street, Cambridge, MA 02138 USA}
\email{yusheng@math.harvard.edu}
\author{Arnau Padrol}
\address{Sorbonne Universit\'e, Institut de Math\'ematiques de Jussieu - Paris Rive Gauche (UMR 7586), Paris, France}
\email{arnau.padrol@imj-prg.fr}

\thanks{Work of the first author  was partially supported by  NSF grants DMS-1100373,  DMS-1401224
and DMS-1701576.  A.\ Padrol thanks the support of the grant ANR-17-CE40-0018 of the French National Research Agency ANR (project CAPPS), as well as the program PEPS Jeunes Chercheur-e-s 2017 of the INSMI (CNRS)}
\date{July 14, 2018, v116c}

\begin{document}
\begin{abstract}
Moser's shadow problem asks to estimate the shadow function~$\shb(n)$,
which is the largest number such that for each bounded convex polyhedron~$P$ 
with $n$ vertices in $3$-space there is some direction~$\bv$ 
(depending on $P$) such that,
when illuminated by parallel light rays from infinity in direction $\bv$,
the polyhedron casts a shadow having at least~$\shb(n)$ vertices.
A general version of the problem allows unbounded polyhedra as well,
and has associated shadow function~$\shu(n)$. 
This paper presents correct order of magnitude
asymptotic bounds on these functions. 
The bounded shadow problem has answer 
$\shb(n) = \Theta \big( \log (n)/ (\log(\log (n))\big).$
The unbounded  shadow problem is shown to have the different asymptotic 
growth rate $\shu(n) = \Theta \big(1\big)$. 
Results on the bounded shadow problem follow from 1989 work of Chazelle, Edelsbrunner and
Guibas on the (bounded) silhouette span number $\sib(n)$, defined analogously but with arbitrary light sources.  We complete the picture by showing that
the unbounded silhouette span number $\siu(n)$ grows as $\Theta \big( \log (n)/ (\log(\log (n))\big)$.
 \end{abstract}

\maketitle

\begin{center}
  {\em  In  memory of Leo Moser (1921\,--\,1970)} 
\end{center}

\section{Introduction}\label{sec:intro}

This paper gives complete answers to 
 several different variants of a problem raised in 1966  in  an influential list of 
problems in discrete and combinatorial geometry made by  Leo Moser~\cite{Mo66}, later reprinted
in 1991 ~\cite{M91}. Moser's life and work are described in  \cite{Mo05}, \cite{Wy72}.

Problem 35 of Moser's list is as follows.\footnote{We have changed
the original notation $f$ to $\shf$ in stating Problems~\ref{main-prob}
and~\ref{main-prob2}.}
\begin{prob}
\label{main-prob}
Estimate the largest $\sh=\sh(n)$ such that every convex polyhedron of $n$ vertices has an orthogonal projection 
onto the plane with $\sh(n)$ vertices on the `outside'.
\end{prob}

A nearly equivalent  problem was  formulated  in a 1968 paper
 of G.  C. Shephard~\cite[Problem VIII]{Sh68b}.
\begin{prob}
\label{main-prob2}
 Find a function $\shf(v)$ such that every convex polyhedron with $v$ vertices possesses a projection which is an
$n$-gon with $n \ge \shf(v)$.
\end{prob}

This problem has been called \defn{Moser's shadow problem} (\cite[p. 140]{CEG89},~\cite[Problem B10]{CFG91}),
because such projections can be viewed as the shadow of the polyhedron
cast by parallel light rays coming 
from a light source ``at infinity." 

The problem can be formulated in  two variants, depending on whether or not unbounded
 polyhedra are allowed.    Shephard's version of the problem ~\cite{Sh68a,Sh68b}  definitely restricts to bounded polyhedra
since  he treats polyhedra that are the convex hull of a finite set of points.  
Following standard terminology  
 such a convex hull is called a \defn{polytope} (\cite[p.4]{Z95}).
 Moser's original problem statement
 does not  explicitly indicate  whether polyhedra are required to  be bounded,
 though he probably had bounded polyhedra in mind.
 In any case the  unbounded version of the problem is of interest 
 because  polyhedra  defined as  intersections of half-spaces naturally arise 
 in linear programming, and certain linear programming algorithms
 have an interpretation in terms of shadows.

In this paper we consider both the bounded and unbounded case.
To distinguish the bounded case from the general (unbounded) case 
 we let~\defn{$\shb(n)$} denote the minimal value over bounded polyhedra (i.e., $3$-polytopes) having
  $n$ vertices, and~\defn{$\shu(n)$} denote the minimal value 
  allowing unbounded polyhedra with $n$ vertices as well (counting only bounded vertices).
We call \defn{Moser's shadow problem} the problem of determining the growth rate of~$\shb(n)$.
We also  formulate in analogy \defn{Moser's unbounded shadow problem}, 
which concerns the growth rate of~$\shu(n)$.

A related problem, the \defn{silhouette span problem}, was formulated by Chazelle, Edelsbrunner and Guibas in 1989~\cite{CEG89}. It  is a variant of the shadow problem
that allows more freedom in the location of the light source from which the shadow is cast.
It considers shadows cast by point light sources at finite distance from the polytope.
The corresponding \defn{bounded silhouette span number}~\defn{$\sib(n)$}, is defined analogously as the shadow number, maximizing over all finite locations of the light source. 
It  is also possible to
define  the \defn{unbounded silhouette span number},~\defn{$\siu(n)$}.
Its formal definition is a little subtle,   and is given in Definition \ref{def:unboundedsilhouette}.

These four functions satisfy the following inequalities, 
\[\begin{tabular}{ccc}
  $\sib(n)$ & $\ge$ & $\siu(n)$\\
   \rotatebox{270}{$\ge\;\; $}& & \rotatebox{270}{$\ge\;\;  $}\\
   $\shb(n)$ & $\ge$ & $\shu(n)$.
  \end{tabular}\]
  The two horizontal inequalities hold because the unbounded numbers minimize over a larger set
  than the bounded numbers, for both the shadow problem and the silhouette span problem. 
  The  vertical inequality between silhouette span numbers and shadow numbers
 holds because silhouettes from light sources that are sufficiently far away in the direction
 of a parallel projection  have
at least as many vertices as 
shadows obtained by that parallel projection (see \cite[pp.174--175]{CEG89} for the bounded case; a similar argument holds for unbounded shadows and silhouettes).

 Chazelle, Edelsbrunner and 
Guibas  \cite[Theorem~4]{CEG89} determined the exact asymptotics of the bounded
silhouette span function $\sib(n)$. 

\begin{theorem}\label{thm:BoundedSilhouette}
{\rm (Chazelle-Edelsbrunner-Guibas)}
The bounded $n$-vertex silhouette span number~$\sib(n)$ for $3$-dimensional convex polytopes satisfies  
\begin{equation*}
                   \sib(n) = \Theta\left(\frac{\log (n)}{\log(\log (n))}\right).
\end{equation*}
\end{theorem}

In this paper, our object is to determine the asymptotic growth rates of the other three functions
$\shb(n)$, $\shu(n)$ and $\siu(n)$, as $n \to \infty$.  In particular, the original Moser shadow problem 
corresponds to   $\shb(n)$.

Our first result puts on  record a complete solution to 
Moser's shadow problem in the bounded polyhedron case.

\begin{theorem}\label{thm:MosersShadowProblem}
The bounded $n$-vertex shadow number~$\shb(n)$ for $3$-dimensional convex polytopes satisfies  
\begin{equation*}
                                                   \shb(n) = \Theta\left(\frac{\log (n)}{\log(\log (n))}\right).
\end{equation*}
\end{theorem}

As we shall  explain below, this result  should  be attributed to  
Chazelle, Edelsbrunner and Guibas, in 
the sense that  all the ingredients for a proof are present in their 1989 paper~\cite{CEG89}.
However, although they mentioned the shadow problem they did not point out  that their
results implied a solution. See Section~\ref{subsec:related} below for more details.

In Section \ref{sec:shadow} 
we present  a complete  self-contained 
proof of Theorem \ref{thm:MosersShadowProblem}.
For the lower bound result,  
 in place of polarity used in  the approach from~\cite{CEG89} sketched in Section~\ref{subsec:related}, 
we use   the spherical image map, introduced in Section~\ref{subsec:GreatCircleSpan},
 together with central projections. This method combines the two first steps 
 in the approach of Chazelle, Edelsbrunner and Guibas. The 
 crucial remaining step in the lower bound, stated below  as Theorem \ref{thm:lblinespan}, is 
  due to \cite{CEG89}, and we supply details 
 for completeness.
  To obtain the upper bound we present a construction that is
  direct and simpler than the one  in \cite{CEG89}, but applies only to
  the shadow problem and not to the silhouette span problem.
  See Section \ref{subsec:ub1} for more details.
 
The remainder of the paper is devoted to 
 the unbounded polyhedron
versions of the shadow and silhouette span problems.

 In Section~\ref{sec:unbounded} we prove 
 that the unbounded shadow function $\shu(n)$ is eventually constant.

\begin{theorem}\label{thm:UnboundedMoserShadowProblem}
The unbounded $n$-vertex shadow number~$\shu(n)$ for $3$-dimensional convex polyhedra satisfies  
\begin{equation*}
                                                   \shu(n) = \Theta(1).
\end{equation*}
In fact $\shu(n) = 3$ for all $n\geq 3$ (and $\shu(1) = 1$ and $\shu(2)=2$).
 \end{theorem}
 
 The upper bound is obtained  by an explicit construction.
 The  lower bounds are obtained with a simple argument
 using spherical images and central projections. 

In Section~\ref{sec:unboundedsilhouette} we treat the unbounded
version of the silhouette span problem.
There is a subtlety  in 
 generalizing the definition of silhouette span to unbounded polyhedra.
 Certain edges visible in an unbounded shadow may not correspond to a  edge of the unbounded
 polyhedron itself.   Our definition, which in the bounded polyhedron case is equivalent to that
 used in \cite[Sect. 5.3]{CEG89},
 allows  as potentially visible  edges corresponding to 
 the recession directions of the unbounded polyhedron. See Definition~\ref{def:unboundedsilhouette}.
We  obtain  the following result, which shows the order of magnitude of the silhouette span number does not decrease when one allows unbounded polyhedra.

\begin{theorem}\label{thm:UnboundedSilhouetteSpan}
The unbounded $n$-vertex silhouette span number~$\siu(n)$ for $3$-dimensional convex polyhedra satisfies  
\begin{equation*}
           \siu(n) =\Theta\left(\frac{\log (n)}{\log(\log (n))}\right).
\end{equation*}
\end{theorem}
 
This result is proved by reduction to the bounded silhouette span case.
Notice that our results show that the shadow problem and silhouette 
span problems have different growth rates
in the unbounded case (in contrast with the bounded case, where both coincide).

\subsection{Related work}\label{subsec:related}

After Moser's original formulation in 1966, the problem was restated several times~\cite{CFG91,Mo66,M91,Sh68b}. The  problem book of Croft, Falconer and Guy~\cite[Problem B10]{CFG91} reports that Moser conjectured $\shb(n)=\sO(\log (n))$
and it sketches the construction of a polytope whose shadow number is of this order of magnitude.
Shephard~\cite[Problem~VIII]{Sh68b} did not conjecture a value for~$\shb(n)$. However, in 
the dual formulation terms of sections~\cite[Problem~VI]{Sh68b}, he proposed a lower bound 
for the silhouette span problem of the form $n^\alpha$ for some constant $0<\alpha<1$.

The 1989  paper 
of Chazelle, Edelsbrunner and Guibas \cite{CEG89} treated a diverse
set of problems concerning the combinatorial and computational complexity of diverse stabbing problems in dimensions two and three. 
In particular the \defn{silhouette span problem}  consists in finding 
the maximal number~$\sib(n)$
such that for each $3$-polytope with~$n$ vertices there is a point from which the silhouette cast has at least~
$\sib(n)$ vertices.  Their approach to the silhouette span problem
 (in the bounded case) exploited the polarity operation (with respect to a point),
 which is a duality operation that interchanges points
 and hyperplanes and preserves incidences (see~\cite[Section~5.1]{Mat2002} for a brief introduction). 
 It  associates a polar polytope~$\pol{P}$ to each polytope~$P$ containing the origin in its interior;
 each point of the polar polytope  corresponds to a particular hyperplane $H$ in $\PP^3$ %
 that lies outside $P$
 in the sense of being disjoint from its interior.
 The polar polytope $\pol{P}$ also contains the origin in its interior; this origin corresponds to
 the plane at infinity in the space of $P$. %
 For each point $p\in \RR^3\smallsetminus P$, its associated plane~$H_p$ is a plane
 that intersects~$\pol{P}$ and does not contain the origin. It is not hard to see 
 that the number of vertices of the silhouette of $P$ as seen from $p$ %
 coincides with the number of edges of the 
 intersection $H_p\cap \pol{P}$. 
 Hence, finding the silhouette span of $P$ is equivalent to finding 
 the maximal number of facets 
  of~$\pol{P}$ which can be intersected with a plane. This 
 problem is referred to in~\cite{CEG89} as the \defn{cross-section span 
 problem}. 
This problem is actually another of the problems in Shephard's list~\cite[Problem~VI]{Sh68b}.
The cross-section span problem is then solved in \cite[Sect. 5.2]{CEG89} using a $2$-dimensional reduction.
The fact that $\sib(n) \ge \shb(n)$, yielding  an upper bound
for $\shb(n)$, was noted on\cite[pp.174--175]{CEG89}.

As we remarked above,  \cite{CEG89}   contains ingredients 
sufficing to prove a lower bound for~$\shb(n)$.
Polarity is actually a kind of projective duality operation. 
Note that shadows from orthogonal projections are the same as silhouettes from points at infinity. The polars of  points at infinity are planes
 through the origin. Hence, the shadow number of~$P$ coincides with the maximal size 
 (number of edges) of a section of~$\pol{P}$ with a plane containing the origin.
  Although in~\cite{CEG89} the authors only claim results for the silhouette span problem and
 the cross-section span problem, their lower bound proof for cross-section span only uses hyperplanes 
 through the origin~\cite[Lemma~5.1]{CEG89}. 
  Therefore, their lower bound
 of order $\Omega({\log (n)}/{\log(\log (n))})$ is also valid for Moser's shadow problem. Thus Theorem \ref{thm:MosersShadowProblem}
follows from the results in \cite{CEG89}.
However, the relevant bound in Lemma 5.1 is stated for an unnamed function $c_d^{\ast}(n)$ and
their paper did not remark on its consequences for the shadow problem, which has been considered open until now.

Very recently Glisse et al. \cite{GlisseLazardMichelPouget2016} studied
the  expected shadow number of a random $3$-polytope obtained by a Poisson point process on the sphere  and showed it to be of order~$\Theta(\sqrt{n})$.

\subsection{Higher-dimensional generalized shadow problems}\label{subsec:higher}
Shadow problems can be generalized to higher dimensions by considering
$k$-dimensional shadows/silhouettes of $d$-dimensional polytopes.  

The special higher-dimensional case of $2$-dimensional projections of $d$-dimensional polyhedra
has been studied  in connection with linear programming algorithms. 
The {\em shadow vertex simplex algorithm} is a parametric version of 
the simplex algorithm in linear programming introduced by Gass and Saaty~\cite{GS55} in 1955.
The analysis of this algorithm leads to the study of $2$-dimensional shadows of $d$-dimensional polyhedra.
A variant of the algorithm  was studied in detail by  Borgwardt~\cite{Bo82a,Bo82,Bo87,Bo99}.
 Later  Spielman and Teng~\cite{ST04} and Kelman and Spielman~\cite{KS06}, 
studied the shadow vertex simplex algorithm
 in connection with average-case analysis of linear programming problems.

For $k$-dimensional shadows we measure size as  the number of vertices visible in 
the shadow; other measures of size may be also considered for $k \ge 3$.
Problems on the size of $k$-dimensional projections of $d$-dimensional polyhedra,
can be translated into problems of intersecting $d$-polyhedra with $k$-dimensional subspaces,
using arguments similar to those given in Section~\ref{sec:definitions} .

Several different types
 of higher-dimensional shadow problems can be considered: worst-case, average-case and minimax case. \begin{enumerate}
\item
{\em Worst case problems} concern the problem of {\em maximizing} shadow numbers for
the a fixed number of vertices. 
The worst case behavior of the shadow vertex method  is related to polyhedra having large shadows,
 For dimension $d=3$ it is easily seen that for all $n \ge 4$  there are 
 polyhedra having all vertices visible in 
a shadow: one may take a suitable oblique cone over a base that is an $(n-1)$-gon. 
Amenta and Ziegler~\cite{AZ98}
 and  G\"{a}rtner, Helbling, Ota and Takahashi~\cite{GHOT13} (see also~\cite{GJM12}) present constructions  of  bad examples of $2$-dimensional shadows in all higher dimensions $d$.   
\item
{\em  Average case problems} concern the average size of $k$-di\-men\-sion\-al shadows
  taken with respect to some measure on the set of directions. 
  Such problems for $2$-dimensional shadows arose from the average case analysis of the
  shadow vertex algorithm. In the 1980's  Borgwardt~\cite{Bo82a, Bo82, Bo87,Bo99} developed
 a polynomial time average case analysis of the variant of the simplex method for linear programming
 that uses the shadow vertex pivot rule.   The shadow vertex simplex algorithm later provided the
  fundamental example used in Spielman and Teng's~\cite{ST04}
  theory of smoothed analysis of
 algorithms. Their analysis requires obtaining some control on the (average) size of shadows,
 as a function of the numbers of variables and constraints in the linear program.
  Further developments of smoothed analysis
 are given in Despande and Spielman~\cite{DS05} and Kelner and Spielman~\cite{KS06}.
 
\item
 {\em Minimax shadow problems} for $2$-dimensional shadows
 in dimensions $d \ge 4$ generalize the shadow problem treated
 in this paper. T\'oth~\cite{T08} has studied line stabbing numbers of convex subdivisions in all dimensions,
extending the analysis of Chazelle et al.~\cite{CEG89}.
His lower bounds induce lower bounds for $2$-dimensional shadow numbers of $d$-polyhedra, however
his examples for  upper bounds are not face-to-face, and hence do not arise from convex polytopes. 
 \end{enumerate}
 
 The general minimax problem  for $k$-dimensional shadows is:

  \begin{prob}
  Estimate the growth rate of the maximal number $\sh_b(n,d,k)$ (resp. $\si_b(n,d,k)$)
  such that every $d$-polytope with $n$ vertices has 
  a $k$-dimensional shadow (resp. silhouette) with $\sh_b(n,d,k)$ (resp. $\si_b(n,d,k)$)
  vertices. Do the same for maximizing over all $d$-polyhedra $\sh_u(n,d,k)$ (resp. $\si_u(n,d,k)$.)
 \end{prob}

To our knowledge all these minimax problems are open in dimensions $d \ge 4$; and so are the analogue 
 silhouette span questions.
\subsection{Plan of the Paper}

Section~\ref{sec:definitions}
gives definitions and 
relates Moser's bounded and unbounded shadow problems to
stabbing problems for spherical and Euclidean polyhedral
subdivisions. These reductions are used in the subsequent sections.
Section~\ref{sec:shadow} contains a full proof
for Theorem~\ref{thm:MosersShadowProblem},  Moser's shadow problem for bounded polytopes. 
The unbounded case of the shadow problem  is treated in Section~\ref{sec:unbounded}, which gives a proof of Theorem~\ref{thm:UnboundedMoserShadowProblem}. Section \ref{sec:unboundedsilhouette}
formulates and treats the unbounded case of the silhouette span problem. 

\section{ Shadows, silhouettes, great circles and stabbing lines}\label{sec:definitions}

We follow the terminology for convex polytopes in Ziegler~\cite[pp.~4--5]{Z95}, 
and define a \defn{polyhedron} in $\R^d$ to be a finite intersection of closed half-spaces, which may be unbounded,
and a \defn{polytope} in $\R^d$ to be the convex hull of a finite set of points; that is, a bounded polyhedron.
Faces of dimensions $0$, $1$ and $d-1$ of a $d$-dimensional polyhedron are called \defn{vertices}, \defn{edges},
and \defn{facets}, respectively.
We say that a polyhedron is \defn{pointed} if it does not contain a full line.
This paper exclusively considers the $3$-dimensional case $\R^3$. 
\subsection{Shadow numbers and silhouette span numbers}\label{subsec:Shad-Silh}

We first define shadows in terms of parallel projections  in a given direction.
\begin{defi}\label{def:shadow-def}
A \defn{shadow} of a (possibly unbounded) polyhedron~$P$ in $\R^3$ is
the image of~$P$ under a %
(possibly oblique) affine projection $\pi_V: \RR^3 \to V$ onto
a two-dimensional affine flat $V$.
The \defn{shadow number} $\sh(P)$ of $P$ is the maximum number of vertices on the boundary of one of its shadows.
\end{defi}

The $2$-dimensional affine subspace $V$ that is the range of the projection $\pi_V$ serves as a ``screen"
on which the shadow $\pi_V(P)$ appears; it is in general a (possibly unbounded) polygon.
In this definition we may  restrict $\pi_V$ to be  
orthogonal projections 
onto a linear subspace $V$ perpendicular to a given $\bv\in\SSS^2$, which we define to be the 
\defn{shadow} %
in direction~$\bv$. 
Therefore, the shadow number $\sh(P)$ of $P$ is the maximal number of vertices visible in shadows of~$P$ obtained
by parallel projection in any direction $\bv$
$$
\sh(P) := \max \{ \sh(P; \bv): \, \bv \in \SSS^2\}.
$$
This definition of shadow number $\sh(P)$ makes sense for both bounded and unbounded polyhedra.

Alternatively, the {shadow number} $\sh(P)$ of~$P$ can also be interpreted as the 
maximal number of $1$-dimensional faces of
the ``cylinder'' resulting from the Minkowski sum $P+\RR \bv$, varying over all 
directions $\bv$. 

We now define the bounded shadow number function as a min-max quantity.

\begin{defi}\label{def:bounded-shadow-num}

The \defn{$n$-vertex bounded shadow number} $\shb(n)$
is given by 
$$
\quad \,\,\, \shb(n) : = \min \{ \sh(P):\mbox{ $P$ is a bounded $3$-polyhedron with $n$ vertices} \}.
$$ 
The \defn{$n$-vertex unbounded shadow number} $\shu(n)$ 
is given by
$$ 
\shu(n) : = \min \{ \sh(P):\mbox{ $P$ is a $3$-polyhedron with $n$ (bounded) vertices} \} \quad\quad;\\
$$
\end{defi}
Note that~$\shu(n)$ could be  referred to as the {$n$-vertex (general) shadow number}, because it 
contemplates bounded and unbounded polyhedra, but we chose this notation to highlight the contrast
with the bounded case. Of course, $\shu(n)\leq\shb(n)$.

\begin{figure}[htpb]
\centering
\includegraphics[width=.75\linewidth]{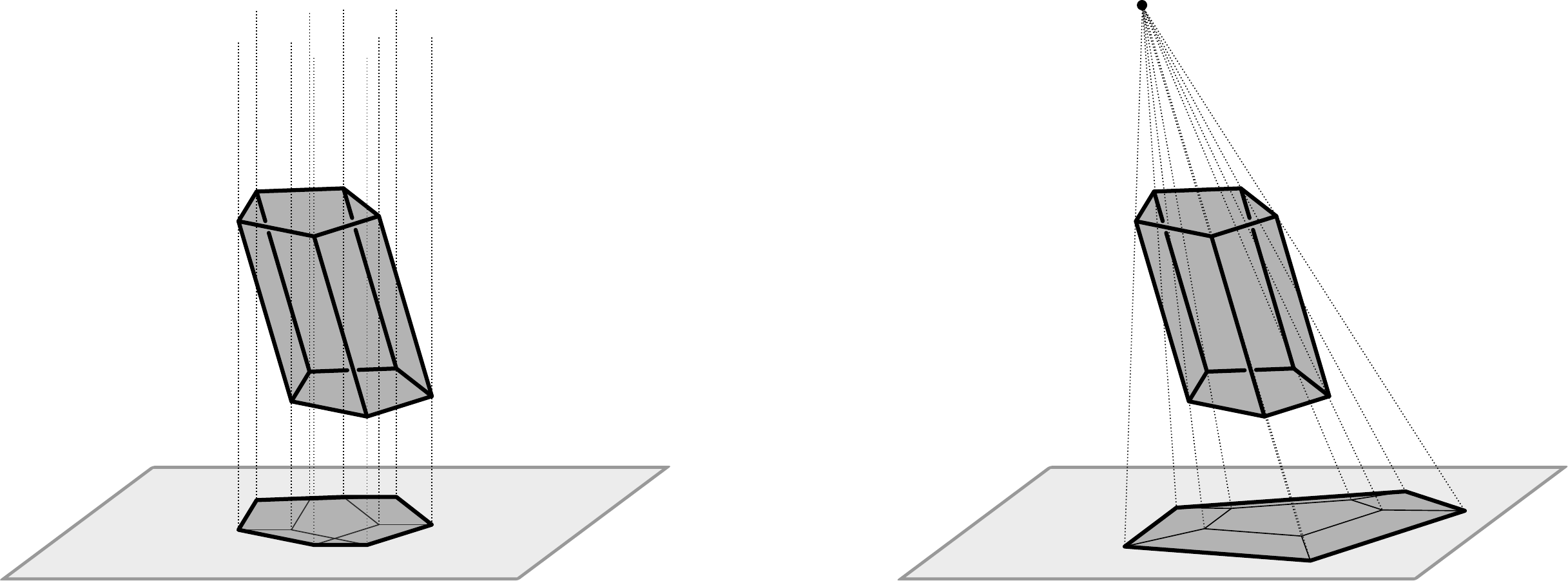}
\caption{A shadow and a silhouette of a polytope.}
\label{fig:shadowVSsilhouette}
\end{figure}

Now we turn to silhouette span. The definition of silhouette span of a bounded polyhedron $P$ given in~\cite[Section~5.3, p.~174]{CEG89},
is an intrinsic definition as a subset of the boundary of $P$. Here we use an alternative definition, equivalent as far as the bounded silhouette span is concerned, that parallels the ``cylinder''
definition of shadow numbers and is better suited for unbounded polyhedra. 

\begin{defi}\label{def:unboundedsilhouette}
Let $P\subset\R^3$ be a (possibly unbounded) polyhedron and $p\in \R^3$ a point outside $P$, and let 
$$C_p(P)=%
\overline{\{p+\lambda \bv\ : \ \bv\in P-p, \lambda \geq 0\}}$$
 be the closure of the
cone with apex $p$ spanned by~$P$. A \defn{silhouette} %
of $P$ with respect to $p$ is a section of $C_p(P)$ with a transversal hyperplane (for example, a hyperplane separating $p$ from~$P$).
The \defn{size} of a silhouette is its number of vertices (in bijection with the rays of the cone), and the \defn{silhouette span}
$\si(P)$ is the size of the largest silhouette of~$P$.
\end{defi}

In \cite{CEG89}, they define the silhouette of a bounded polytope~$P$ 
with respect to a point~$p$ outside~$P$ as the collection of faces~$F$ of~$P$ that allow a supporting
plane~$H$ of~$P$ such that~$p$ lies in~$H$ and~$F$ is in the relative interior of~$P \cap F$; and measure its size as its number of vertices.
To avoid confusion, we may call this the \defn{pre-silhouette} of~$P$ 
with respect to~$p$ (such complexes are sometimes referred to as the \defn{shadow-boundary} of~$P$ from~$p$, see for example~\cite{Sh72}).
When~$p$ is not coplanar with any facet of~$P$, the pre-silhouette is a collection of edges and vertices in the boundary of $P$ (but otherwise it might 
also contains facets). In this case, central projection from $p$ maps the pre-silhouette bijectively to the boundary of the silhouette. Since silhouettes of maximal size are 
always attained from points in general position, both definition give exactly the same silhouette spans.

\begin{figure}[htpb]
\centering
\includegraphics[width=.35\linewidth]{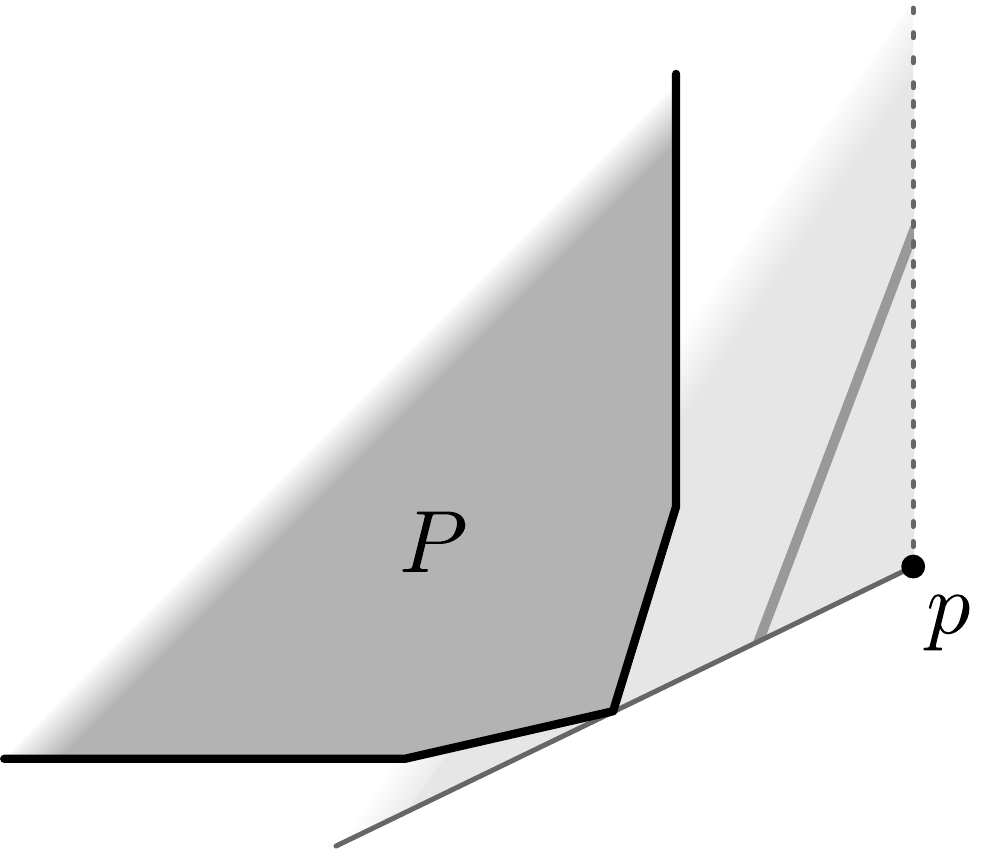}
\caption{A $2$-dimensional unbounded polyhedron~$P$ as seen from a point~$p$. The cone with apex $p$ spanned by~$P$ is not closed, one boundary edge (dotted lines) is missing. Any transversal section of the closure of this cone gives a silhouette, one example is the highlighted segment.}
\label{fig:unboundedsilhouette}
\end{figure}

However, this definition of pre-silhouettes is not well adapted for unbounded polyhedra. 
If $P$ is unbounded, we wish to consider 
also as part of  the silhouette those faces of the recession cone that are visible from~$p$ at infinity.
Indeed, silhouettes can be interpreted by projecting onto a canvas that separates $P$ from a viewer placed at~$p$. Unbounded facets are seen as half-open polytopes, in which part of the boundary may
be  missing, as it corresponds to limit directions at infinity.
An example with missing boundary is sketched in Figure~\ref{fig:unboundedsilhouette} for the planar case.

Our definition includes this extra boundary (this is why the closure is needed in the definition of~$C_p(P)$).
In order to reformulate the definition of pre-silhouettes to this set-up, one should consider also some extra unbounded edges of $P$
in the directions of the recession cone. To each such 
unbounded edge  it adds  a ``vertex at infinity" for
such edges.  The \defn{silhouette size} of $P$ viewed from $p$ 
would now count the additional ``vertices at infinity"
included this way.   

We  now define the  bounded  silhouette span function as a min-max quantity.

\begin{defi}
The \defn{$n$-vertex bounded silhouette span number} $\sib(n)$
is given by 
  \begin{align*}
 \sib(n) & : = \min \{ \si(P):\mbox{ $P$ is a bounded $3$-polyhedron with $n$ vertices} \}.
 \end{align*}

 The \defn{$n$-vertex unbounded silhouette span number} $\siu(n)$
is given by 
  \begin{align*}
 \siu(n) & : = \min \{ \si(P):\mbox{ $P$ is a is a $3$-polyhedron with $n$ (bounded) vertices} \}.
 \end{align*}
 
 \end{defi}

\begin{rmk}\label{rmk:size}
Definitions~\ref{def:shadow-def} and~\ref{def:unboundedsilhouette} measure the \defn{size} of a polyhedron $P$ in terms of its 
number of vertices, as in Moser's
version of the problem.
Alternatively one can  measure the combinatorial size of $3$-polyhedra in terms of
vertices $v(P)$, edges $e(P)$ or facets $f(P)$, or some combination of
all three terms.
For bounded polyhedra, all three of these complexity measures are related
within a linear factor by Euler's formula
(cf.~\cite[pp.~189--190]{Gruenbaum}).
Hence, the bounded shadow number and bounded silhouette span number have the same asymptotic behavior regardless of whether we use  the number of edges, faces, or vertices of $P$ as a measure of its size.

For unbounded polyhedra, there is no lower bound
relating $v(P)$ to $e(P)$; there can be one vertex and arbitrarily many edges.
However the construction in the proof of Theorem~\ref{thm:UnboundedMoserShadowProblem} has an unbounded number of vertices, edges, and facets; 
and Theorem~\ref{thm:UnboundedSilhouetteSpan} can be easily adapted to give the same asymptotics for $n=e(P)$ and $n=f(P)$, as the missing ``vertices at infinity'' are also visible in the silhouettes.
 \end{rmk}

\begin{rmk}\label{rmk:proj_transformations}
We make some  additional remarks concerning the effect of Euclidean and projective transformations
on shadow number:
\begin{enumerate}
 \item\label{it:norm_eq} The shadow number is a Euclidean invariant of $3$-polyhedra, i.e.\ 
two congruent polyhedra have equal shadow numbers.
In addition, two normally equivalent $3$-polyhedra, i.e.\ which have
the same combinatorial type and identical normal directions to each corresponding face,
have identical shadow numbers (but they are not
necessarily congruent).
 
\item 
The  shadow number of a polytope is not a projective invariant, i.e.\ one
can exhibit  examples of polytopes that are equivalent under  
 a projective transformation in the sense of~\cite[Appendix 2.6]{Z95} which
 have different shadow numbers. 

\item
The  silhouette span number $\sh(P)$ of a $3$-polyhedron $P$ is preserved by
those  projective transformations whose hyperplane at infinity does
not intersect $P$.
 \end{enumerate}
\end{rmk}

\subsection{Reduction via spherical image to great circle span problems}\label{subsec:GreatCircleSpan}

We reduce the shadow problem
in both the bounded and unbounded cases
to (special cases of) a 
dual problem about convex geodesic subdivisions on the standard sphere 
$$
\SSS^2:= \{ u= (u_1, u_2, u_3) \in \R^2: u_1^2 +u_2^2 + u_3^2 =1\}.
$$

We define a \defn{(convex) spherical polyhedron} as a finite intersection of closed 
hemispheres in~$\SSS^{d-1}$, where a \defn{closed hemisphere} is the intersection of~$\SSS^{d-1}$ with a closed halfspace containing the origin in its boundary. We use the term \defn{spherical polygons} to denote 
two-dimensional spherical polyhedra. Their boundary is formed by segments of great circles, called \defn{edges}.

\begin{defi}\label{def:spherical_convex_subdivisions}
Let $U$ be either $\SSS^{d}$ or a convex spherical $d$-polyhedron.
A \defn{spherical polyhedral subdivision} (or \defn{subdivision} for short) 
of $U$ is a finite set of convex spherical $d$-polyhedra (called \defn{regions}),
whose union is $U$ and such that the intersection of any two is a common face.
\end{defi}

Although we could have relaxed the definition of subdivision
without imposing the condition of being face-to-face (as in~\cite{CEG89}), we will be mainly concerned with spherical 
polyhedral subdivisions of $\SSS^{2}$ arising from polyhedra, which are always face-to-face,
and imposing this condition simplifies the exposition.

Let $P\subset\RR^d$ be a polyhedron. Recall that a hyperplane $H$ is called \defn{supporting} for a
face~$F$ of $P$ if $H\cap P=F$ and $P$ is 
completely contained in the closed halfspace opposite to its
\defn{outward unit normal}.
The \defn{spherical image} $\sig(F) \subset \SSS^{d-1}$ of $F$
is the subset of $\SSS^{d-1}$ consisting of all outward unit normal directions to its supporting hyperplanes. 
The union of all these cells is the \defn{spherical image} $\sig(P) \subset \SSS^{d-1}$ of $P$. If~$P$ is 
bounded then $\sig(P) = \SSS^{d-1}$, while if $P$ is pointed but unbounded then $\sig(P)$
is a convex spherical polyhedron contained in some hemisphere of $\SSS^{d-1}$ (see Alexandrov~\cite[Section~1.5, esp. Theorem 3]{A05}).
If $P$ is not pointed, then $\sig(P)$ is completely contained in an equator of $\SSS^{d-1}$ orthogonal to
its linearity space.

The spherical polyhedral subdivision of $\sig(P)$  induced by the spherical images of the faces
of $P$ is the \defn{spherical image subdivision} $\sDD(P)$ of $P$. 
It is the intersection of the normal fan of $P$ with the unit sphere (see~\cite[Section~7]{Z95}).
In particular, normally equivalent polyhedra have 
the same spherical image subdivision, cf. Remark~\ref{rmk:proj_transformations}\eqref{it:norm_eq}.
However, in general, polyhedra having the same combinatorial type 
will have different spherical images.

Hence, the spherical image subdivision of a (pointed) $3$-polyhedron~$P$ is either~$\SSS^2$ 
or a convex spherical polygon; and the spherical images of its facets, edges, and vertices are respectively 
points, segments of a great circles, and spherical polygons.
\begin{rmk}\label{rmk:polyhedra_vs_scs}
Spherical image subdivisions $\sDD(P)$ of bounded $3$-polyhedra
are spherical polyhedral subdivisions of~$\SSS^2$.
However, the reciprocal does not hold.
Indeed, the spherical image subdivision of a polytope is always a regular subdivision of
a vector configuration and there are subdivisions of vector configurations
that are not regular. See Sections 9.5 and 2.5 of de Loera et al.~\cite{dLRS10} and
also Connelly and Henderson~\cite{CH80}.
\end{rmk}

Our interest on spherical subdivisions is motivated
by the fact that it is possible to read off the shadow number of a polytope
from its spherical image subdivision.

\begin{defi}\label{def:greatcirclespan}
Let $U$ be either~$\SSS^2$ or a convex spherical polygon, 
and let~$\sD$ be a spherical polyhedral subdivision of $U$.
For each great circle $C$ in~$\SSS^2$, the intersection
$C\cap\sD$ induces a spherical subdivision of $C\cap U$.
The \defn{great circle span} $\cs(\sD)$ of $\sD$ is the maximal number 
of regions of a subdivision $C\cap\sD$ obtained this way.
 \end{defi}
Although one is tempted to define $\cs(\sD)$ simply as the maximal number
of regions whose interiors are intersected by a great circle, it is important to
also take into account the cases where a great circle goes along an edge. Otherwise
the following lemma would not hold in some degenerate cases (cf. Remark~\ref{rmk:degeneratecase}).

\begin{lem}\label{lem:equivalencegreatcircleshadow}
Let $P$ be a pointed
polyhedron in $\R^3$, and let~$\sDD(P)$ be its induced spherical image subdivision
of $U = \sigma(P)$.
Then the shadow number of $P$ coincides with the great circle span of $\sDD(P)$:
\begin{equation*}
\sha(P) = \cs(\sDD(P)).
\end{equation*}
\end{lem}

\begin{proof}
For $\bv\in \SSS^2$, let $C_\bv\subset\SSS^2$ denote the great circle perpendicular to $\bv$,
and $\pi_\bv$ the orthogonal projection along $\bv$. We will show that
the number of vertices of $\pi_\bv(P)$ coincides with the number of regions
of $\sDD(P)\cap C_\bv$. This follows essentially from~\cite[Lemma~7.11]{Z95},
which shows that the spherical image subdivision of $\pi_\bv(P)$ coincides with~$\sDD(P)\cap C_\bv$.
Therefore, the maximal number of vertices of a projection coincides with the maximal number of arcs of the 
subdivision induced on a great circle.

Indeed, the arcs of $\sDD(\pi_\bv(P))$ correspond to the sets of outer normal vectors of supporting
hyperplanes for each of the vertices of $\pi_\bv(P)$. Notice that, if $v$ is a vertex of $\pi_\bv(P)$, and $F=\pi_\bv^{-1}(v)$ its pre-image,
then the pre-image of each supporting hyperplane for $v$ in~$\pi_\bv(P)$ is a supporting hyperplane
for $F$ in $P$ whose normal vector is orthogonal to $\bv$. Hence, $C_\bv$ intersects the spherical image of~$F$
in the segment of great circle corresponding to $\sig(\pi_\bv(P))$. 
This argument is reversible. Each segment of great circle of $C_\bv\cap \sDD(P)$ arises from the intersection 
of $C_\bv$ with the spherical image of a face~$F$. The supporting hyperplanes corresponding to these intersection
points are orthogonal to $\bv$, which implies that they are the pre-images of supporting hyperplanes for $\pi_\bv(F)$ on $\pi_\bv(P)$.
\end{proof}

\subsection{Reduction by central projection to stabbing number problems}\label{subsec:Stabbing}

We relate great circle span problems on spherical polyhedral subdivisions of~$\SSS^2$
to  a family of  line span (stabbing number) problems on Euclidean polyhedral subdivisions of $\R^2$
using central projection.

\begin{defi}\label{def:convexsubdivision}
Let $P$ be either $\RR^{d}$ or a $d$-polyhedron.
A \defn{Euclidean polyhedral subdivision} (or \defn{subdivision} for short) 
of $P$ is a finite set of $d$-polyhedra (called \defn{regions}),
whose union is~$P$ and such that the intersection of any two is a common face.
\end{defi}

\begin{defi} \label{def:linespan}
Let $P$ be either $\RR^2$ or a convex polygon, 
and let~$\sE$ be a Euclidean polyhedral subdivision of $P$.
The \defn{line span} or \defn{stabbing number} $\ls(\sE)$ of $\sE$ is 
the maximal number of segments of the restriction of $\sE$ to a line.
 \end{defi}

\defn{Central  projection}, also called \defn{gnomonic projection},
maps an open hemisphere of a sphere~$\SSS^2$ bijectively to a plane. 
Let~$\SSS^2$ be the standard sphere in $\R^3$, 
and let $\SSS_{-}^2 := \SSS^2 \cap \{(x_1,x_2,x_3):x_3<0\}$ denote
 the lower open hemisphere in the last coordinate. 
 Let~$H$ be the plane $\{(x_1,x_2,x_3):x_3=-1\}$. We define the \defn{central 
 projection} $\gamma:\SSS_{-}^2 \longrightarrow H$ by mapping $\bv\in \SSS_{-}^2$ 
 to the unique intersection point $u\in H$ of the line through $0$ and $\bv$ with~$H$.

Now let $\sD$ be a spherical polyhedral subdivision of~$\SSS^2$. 
Let $\sD_{-}$ the convex subdivision of the open hemisphere $\SSS_{-}^2$, 
obtained by intersecting each region of $\sD$ with $\SSS_{-}^2$. Its image under  the gnomonic projection
 $\gamma:\SSS_{-}^2\longrightarrow H$
 is then a  (Euclidean) polyhedral subdivision $\gamma(\sD_{-})$ of~$H$.

This map yields the following relationship between great circle spans and 
line spans, cf.~\cite[Lemma 5.1]{CEG89}.

\begin{lem}\label{lem:equivcirclelinespan}
Let $\sD$ be a spherical subdivision of~$\SSS^2$ or of a convex spherical polygon~$U$, 
with $\sD_{-}$ and $\sD_{+}$ being its restrictions to the upper and lower open hemispheres, 
and let the Euclidean polyhedral subdivsions $\sE_-=\gamma(\sD_-)$ and $\sE_+=\gamma(-\sD_+)$ be
their respective central projections.
Then, 
\begin{equation}\label{eq:ub_csls}
\cs(\sD) \geq \ls(\sE_-).
\end{equation}
If moreover the equator does not contain any edge of $\sD$, then 
\begin{equation}\label{eq:lb_csls}
\cs(\sD) \leq \ls(\sE_-)+\ls(\sE_+).
\end{equation}
Finally, if $U\subset\SSS_-^2$, then
\begin{equation}\label{eq:eq_csls}
\cs(\sD) = \ls(\sE_-).
\end{equation}
\end{lem}

\begin{proof}
The proof  is immediate by noting that 
the central projection of (the restrictions to $\SSS^2_ -$ of) great circles 
and convex spherical polygons are respectively lines and Euclidean polygons (when not empty).
\end{proof}

\begin{rmk}\label{rmk:degeneratecase}
 A subtlety  in Lemma \ref{lem:equivcirclelinespan}
 lies in the requirement that the equator should not contain any edge, which is needed
 for the validity of~\eqref{eq:lb_csls}.
 This condition is necessary to cover the cases where the equator might be the 
 single great circle achieving the maximum span. This can only happen in the 
degenerate case of those~$\sD$ whose support's boundary contains a segment of the equator
$\SSS^2\cap\{x_3= 0\}$, %
because for $\sD$ whose support is $\SSS^2$ or is contained in an open hemisphere, 
the great circle span is easily shown to be attained by great circles
in general position (not passing through any vertex), which can be perturbed without affecting the great circle span.

This is not the case for $\sD$ being the closed lower hemisphere, where the equator might be the single
great circle achieving the maximal circle span. For example, consider its subdivision 
whose vertices are the south pole and $2n$ equi-spaced points on the equator, and whose
cells are the triangles joining segments in the equator with the south pole. Its great circle span
is $2n$, whereas any great circle other than the equator intersects at most $n+1$ regions.
This example corresponds to the spherical image subdivision arising from the Cartesian product of
a regular $2n$-gon with a half-line. Only the projection along the direction of the half-line provides
a shadow with $2n$ vertices.
\end{rmk}

\section{Moser's (bounded) shadow problem}\label{sec:shadow}

 This section  presents a self-contained proof for Moser's shadow problem
 for bounded polyhedra. The proof is divided in two  parts,  the lower
 bound and upper bound, which may be read independently.
 
 The solution to the bounded version of 
Moser's shadow problem
uses a lower bound result for \defn{minimal line span} (i.e.\ a minimax stabbing number) proved
by Chazelle, Edelsbrunner and Guibas~\cite[Lemma 3.2]{CEG89}. 
For completeness, we
state and prove this result  below as Theorem \ref{thm:lblinespan}
following closely the original proof. 
It  uses an iterated topological sweep and  captures a
crucial tradeoff explaining why ${\log (n)}/{\log(\log (n))}$ is a lower bound
for the minimal line span. 
 
 The upper bound $$\shb(n) = \sO\left(\frac{\log(n)}{\log(\log(n))}\right)$$
already  follows  from the upper bound for the silhouette span
problem. In~\cite[Lemma~5.15]{CEG89} Chazelle et al.\ construct (the polar dual of a) polytope with~$n$ vertices whose silhouette
from each point of view has size at most $\sO(\log (n)/\log(\log (n)))$. Since shadows can be regarded as a special  kind of silhouettes, and this upper bound matches the lower bound in Theorem~\ref{thm:LB}, this finishes a proof for Theorem~\ref{thm:MosersShadowProblem}.

However the construction in~\cite[Section~5.2]{CEG89} providing an upper bound for the silhouette span problem is very involved, requiring some quite technical steps,
 and is formulated in a polar dual form.  
Constructing lower bound examples for the shadow number  problem is actually simpler.
We  present a  direct construction that establishes the upper bound.

\subsection{The lower bound: minimal line span bound}\label{subsec:lb1}

We aim to prove the following result. 
\begin{theorem}\label{thm:LB} 
The bounded $n$-vertex shadow number~$\shb(n)$ for $3$-dimensional polytopes
satisfies  
\begin{equation*}
 \shb(n) = \Omega\left(\frac{\log (n)}{\log(\log (n))}\right).
\end{equation*}
\end{theorem}

The crucial ingredient is a lower bound for the  line span for convex polygonal subdivisions
due to Chazelle, Edelsbrunner, Guibas \cite[Lemma 3.2]{CEG89} . We supply a version of their proof,
clarifying some points, for the reader's convenience.
\begin{theorem} \label{thm:lblinespan} {\rm (Chazelle, Edelsbrunner, Guibas (1989))} 
Let $\ls(n)$ be the minimal line span of a convex polygonal subdivision of the plane $\R^2$ into $n$ regions. Then 
for all sufficiently large $n$,
\begin{equation*}
\ls(n)  \ge  \frac{\log (2n)}{\log(\log (2n))}.
\end{equation*}
\end{theorem}
 
\begin{proof}
Let $\sE$ be a plane subdivision with $n$ regions.
We fix a coordinate system 
such that no pair of vertices of $\sE$ share the same $y$-coordinate. We claim that,
for each integer $k> 1$, either there is a non-horizontal line~$\ell$ such that the induced subdivision $\ell\cap \sE$ has more than $k$ regions, or
there is a horizontal line $h$ whose induced subdivision $h\cap \sE$ has at least $\log_{6k}(2n)$
regions. Taking $k=\left\lceil\frac{\log(2n)}{\log(\log(2n))}\right\rceil$, this shows
that for every subdivision there is a line stabbing at least
$$\min\left(\left\lceil\frac{\log(2n)}{\log(\log(2n))}\right\rceil, \frac{\log(2n)}{\log\left({\frac{6\log(2n)}{\log(\log(2n))}}\right)}\right)$$
many regions of $\sE$. Notice that, for sufficiently large $n$, 
$$\frac{\log(2n)}{\log\left({\frac{6\log(2n)}{\log(\log(2n))}}\right)}> \frac{\log(2n)}{\log(\log(2n))};$$
which implies the announced bound on the line span.

To prove the claim, assume that no non-horizontal line stabs more than~$k$ regions of $\sE$. 
To find the desired horizontal line, we start with an open axis-parallel rectangle $K_0$ that intersects all the $n$ regions of $\sE$.
In our recursive proof,~$K_i$ will be an open horizontal strip  
bounded on the left and right by convex polygonal paths that we call $L_i$ and $R_i$ (which might share an endpoint but are otherwise disjoint), 
and above and below by horizontal segments (which might collapse in a point)
with the property that there is a non-horizontal open segment~$\ell_i$ that connects a topmost point in the upper edge and 
a bottommost point in the lower edge that lies entirely inside the strip. The intersection with $\sE$ induces a subdivision
$\sE_i$ of~$K_i$, and $n_i$ will be a lower bound for its number of regions.
In $K_0$, the paths $L_0$ and $R_0$ are just 
the left and right edges of the rectangle, any vertical line between the two can play the role of $\ell_0$, and $n_0=n$.

\begin{figure}[htpb]
\centering
\includegraphics[width=.5\linewidth]{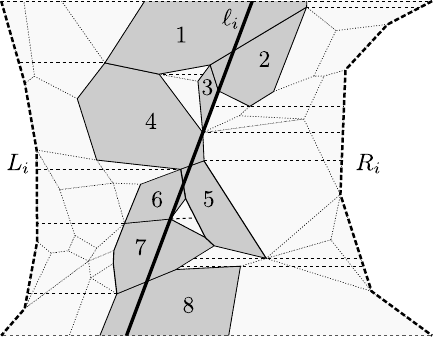}
\caption{A strip~$K_i$ bounded above and below by horizontal lines and left and right by convex polygonal paths $L_i$ and~$R_i$, which are separated by the line $\ell_i$. The regions in~$V_i$ are darker, and numbered according to their order of intersection with $\ell_i$. Several degeneracies are present. For example, the third region shares an edge with $\ell_i$; the lower horizontal segment of the strip bounded left by $4$ and right by $3$ is collapsed into a point; and the strip bounded left by $7$ and right by $5$ consists of two connected components. The strip intersecting more regions of $\sE_i\setminus V_i$ is the one bounded by $L_i$ and region~$7$, and it will be~$K_{i+1}$.   }
\label{fig7}
\end{figure}

For each~$i$, let~$V_i$ be the set of those regions of $\sE_i$ whose interior is intersected by $\ell_i$, or that have an edge on $\ell_i$ and lie at its left.  
Note that the size of~$V_i$ is at most~$k$, because these regions are in bijection with the segments of~$\sE_i\cap\ell_i$.
We now subdivide the complement of~$V_i$ in~$K_i$ into horizontal strips; one of which will be defined as $K_{i+1}$. First, through the topmost
and bottommost vertex of each region in~$V_i$, draw
the longest horizontal segment that does not intersect a region in~$V_i$. This decomposes the complement of~$V_i$ in~$K_i$ 
depending on the first regions in~$V_i$ (or $L_i$ or $R_i$) that are hit when moving horizontally to the left 
and to the right. Note that by construction these strips are bounded at each
side by a convex polygonal path (a part of the boundary of one of the regions in~$V_i$, or $L_i$ or $R_i$). 
If these two boundary paths intersect in an interior edge (which is unique by convexity), the we further subdivide the strip into its two connected components.
In each case, there is an open segment lying entirely inside the region that connects a topmost point and 
a bottommost point. This fact follows from the existence of a separating line
between the two convex polygonal paths.\footnote{This is the crucial point 
in the argument where convexity of all pieces is used.
It no longer holds in the unbounded 
polyhedra  case treated in Section~\ref{sec:unbounded}.}

This subdivision of $K_i\setminus V_i$ has at most $6k$ strips. 
The strips obtained before splitting into connected components can be counted by sweeping a horizontal line
from the bottom: At the beginning it intersects two strips, each time it encounters 
a bottom vertex it enters (at most) two new strips, and one new strip each time it encounters a top vertex. Since each strip can be further subdivided into its 
two connected components, this gives the bound of~$6k$.
There are at most $k$ regions in~$V_i$, and at least $n_i$ regions in $\sE_i$, which means that there are at least $n-k$ regions in $\sE_i\setminus V_i$, each of which intersects at least one strip.
Hence, one of these horizontal strips intersects at least $n_{i+1}:=\frac{n_i-k}{6k}$ regions of $\sE$. We define it to be~$K_{i+1}$
and set $\ell_{i+1}$ to be any of the segments that joins the upper and lower edges.
Notice that 
\[n_i=\frac{n}{(6k)^i} -\sum_{j=1}^i\frac{k}{(6k)^j}=\frac{n}{(6k)^i} -\frac{k-\frac{k}{(6k)^i}}{6k-1}>\frac{n}{(6k)^i} -\frac{k}{6k-1}>\frac{n}{(6k)^i} -\frac{1}{2}.\]
Hence, whenever $i<\log_{6k}(2n)$, we have $n_i> 0$.

Notice how, every horizontal line that goes through the interior of a region in~$V_i$ also intersects the interior of a region in $V_{i-1}$ (the one containing its intersection with $\ell_{i-1}$).
Therefore,
a horizontal line through the interior of a region in $V_{\log_{6k}(2n)-1}$ stabs at least $\log_{6k}(2n)$ distinct regions of~$\sE$. The theorem follows.
\end{proof}

\begin{rmk} \label{rmk:convexitylinespan}
The  hypothesis of convexity made in the statement of 
Theorem \ref{thm:lblinespan} 
is essential; the conclusion can fail badly otherwise.
Indeed  Section~\ref{subsec:unbounded-shad--upper} below gives a construction exhibiting
 (non-convex) polygonal subdivisions of the plane
having  an arbitrarily large number of bounded regions
while still having a uniformly  bounded line span. These subdivisions  have convex bounded regions, 
plus exactly one non-convex unbounded region.
\end{rmk}

\subsection{The lower bound: completion of proof}\label{subsec:lb2}

Now we are ready to complete the proof of the lower bound.

\begin{proof}[Proof of Theorem~\ref{thm:LB}]
Let $P$ be a bounded polytope in $\RR^3$ with $n$ vertices, and let $\sD=\sDD(P)$ be the induced spherical image subdivision
of~$\SSS^2$, which has $n$ regions. By rotating~$P$ (and hence also the subdivision $\sD$)
if needed, we may assume that 
the lower hemisphere~$\SSS_{-}^2$ intersects at least $\lceil n/2\rceil$ regions of $\sD$, and, hence, 
that the central projection $\sE=\gamma(\sD_{-})$ of $\sD_{-}=\SSS_-^2\cap\sD$ has at least $\lceil n/2\rceil$ regions.

By Theorem~\ref{thm:lblinespan}, the line span of $\sE$ is at least 
\[ 
\ls(\sE) = \Omega\left(\frac{\log (n)}{\log(\log (n))}\right).
\]
Combining 
Lemmas~\ref{lem:equivalencegreatcircleshadow} and \ref{lem:equivcirclelinespan} we have
that
\begin{equation*}
\ls(\sE)\leq \cs(\sD)=\sh(P).
\end{equation*}
Thus we have $\shb(n) = \min_{P} \sh(P) = \Omega\left(\frac{\log (n)}{\log(\log (n))}\right).$
\end{proof}

\subsection{The upper bound: direct construction}\label{subsec:ub1}

We will  prove the following result.

\begin{theorem}\label{thm:UB} 
The bounded $n$-vertex shadow number~$\shb(n)$ for $3$-dimensional polytopes satisfies  
\begin{equation*}
                                                   \shb(n) = \sO\left(\frac{\log (n)}{\log(\log (n))}\right).
\end{equation*}
\end{theorem}

The proof will be based on a construction. We will give it in terms of Euclidean polygonal 
subdivisions. An important point in the proof  is to be able to certify that the ones we use
are  gnomonic projections of a spherical subdivision arising from a $3$-dimensional polytope,
which is the polytope we seek to construct.
\begin{defi}
A Euclidean subdivision $\sE$ with~$n$ regions is \defn{liftable} 
if there is a polytope~$P$ with~$n$ vertices such that the central projection of the restriction of its spherical image subdivision 
to the lower hemisphere
coincides with $\sE$, i.e.\ $\sE=\gamma(\sDD(P)\cap\SSS^2_-)$.
\end{defi}

We will repeatedly use three operations. The first pair are classical, based on Steinitz's $\Delta-Y$ operations,
and correspond to the polytope operations of stacking and truncating; the third is a combination of both
these operations.

\begin{figure}[htpb]
\includegraphics[width=.2\linewidth]{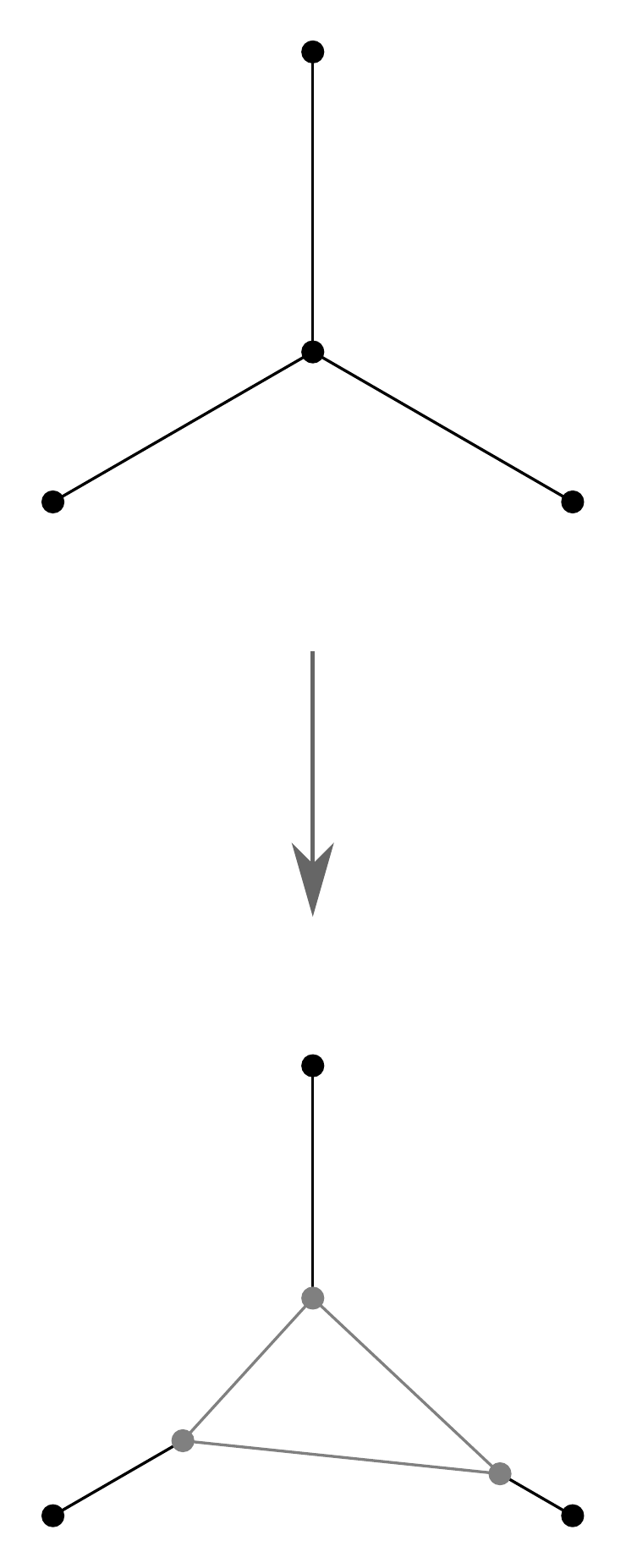}\qquad\qquad
\includegraphics[width=.2\linewidth]{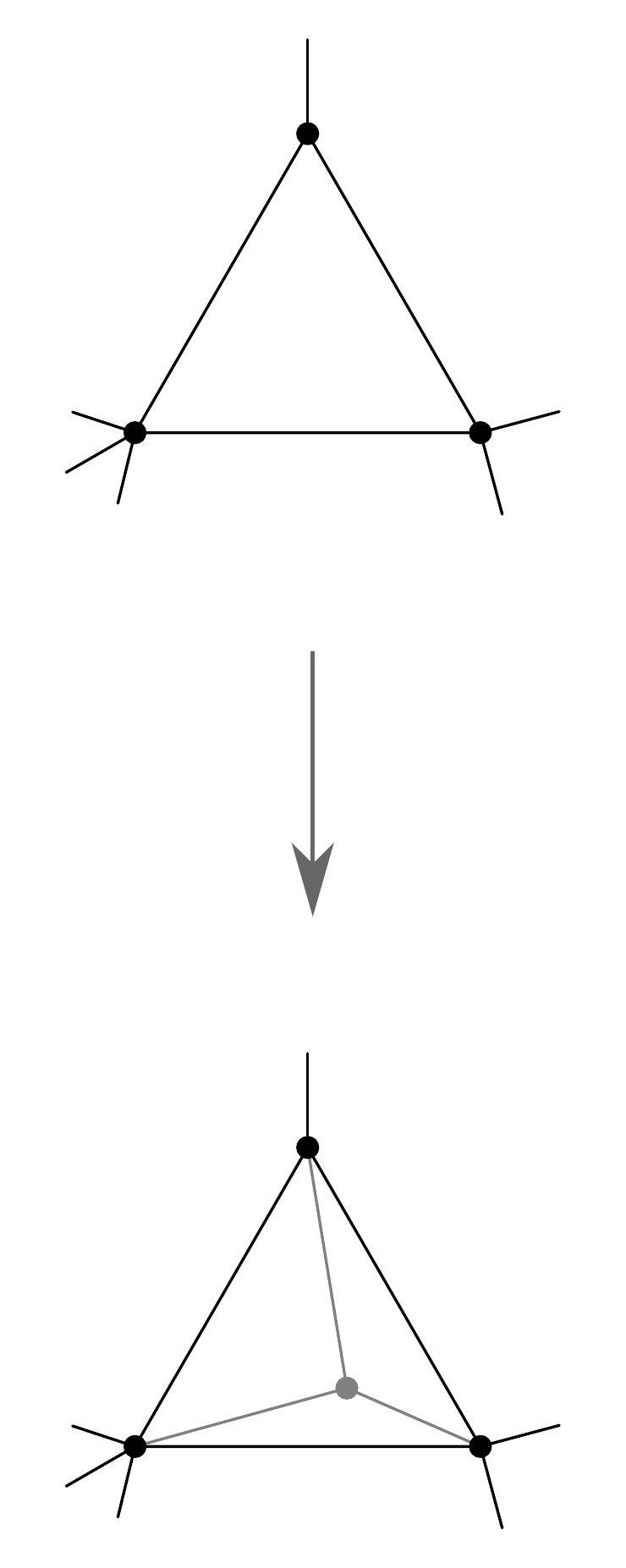}\qquad\qquad
\includegraphics[width=.2\linewidth]{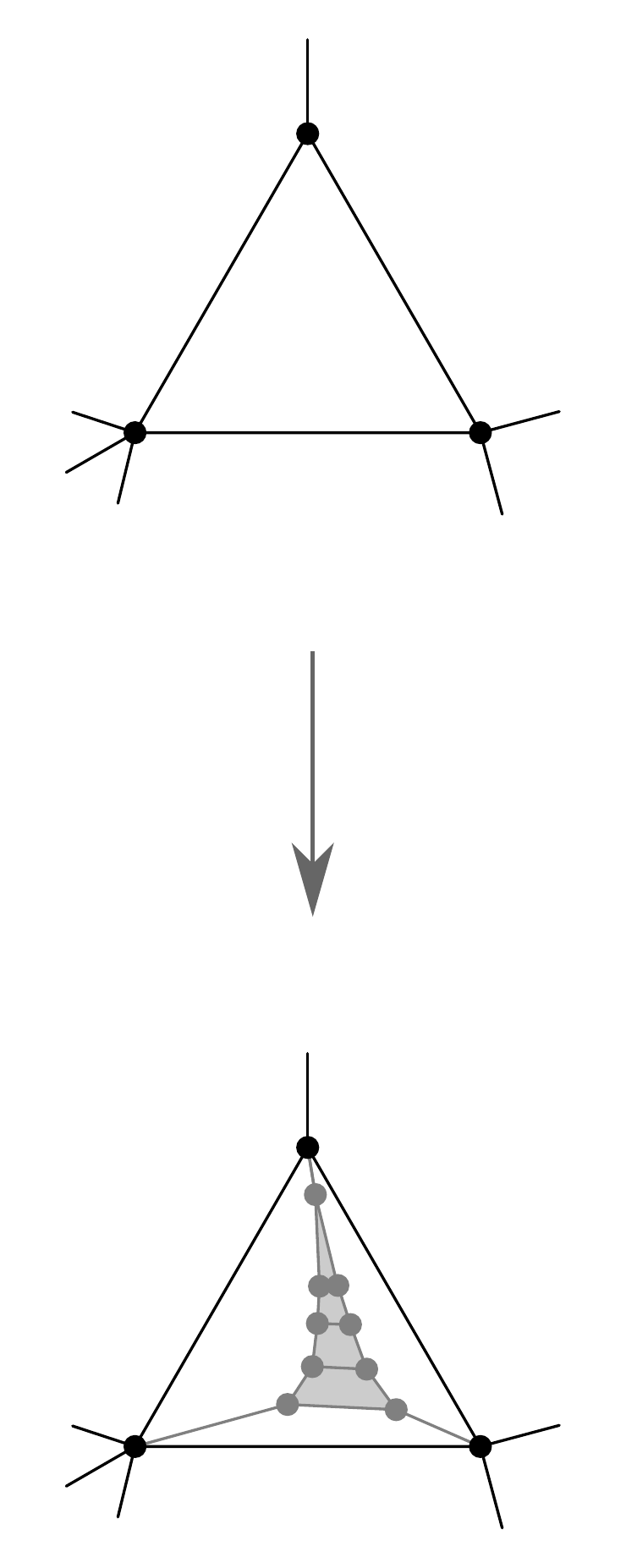}
 \caption{Examples of truncating, stacking and unzipping. The shadowed regions form the spine of the unzipping, which is of length~$4$.}\label{fig:operations}
\end{figure}

\begin{defi}
 Let $\sE$ be polyhedral subdivision of $\R^2$.
 \begin{enumerate}
  \item Let $v$ be a degree-$3$ vertex with neighbors $v_1,v_2,v_3$.
    \defn{Truncating}~$v$ consists in choosing a point~$v_i'$ in the interior of each of the edges $(v,v_i)$
  and 
  adding to~$\sE$ the triangle with 
  vertices $v_1',v_2',v_3'$ (and intersecting the remaining regions with the closure of its complement).
  \item Let $T$ be a triangular region with vertices $v_1,v_2,v_3$. \defn{Stacking} 
  onto $T$ corresponds to adding a vertex $v$ in the interior of $T$ and substituting
  $T$ by the three triangles obtained by joining $v$ with an edge of $T$.
  \item Let $T$ be a triangular region with vertices $v_1,v_2,v_3$.
  \defn{Unzipping} $T$ towards $v_i$ is an operation that consists in first stacking onto $T$ 
  and then successively truncating the newly created vertex that is connected to $v_i$.
  Its \defn{length} is the number of truncations, and the regions created with
  the truncations are the \defn{spine}.
 \end{enumerate}
See Figure~\ref{fig:operations} for an example.
\end{defi}

\begin{lem}\label{lem:liftable}
 Truncating and stacking, and hence also unzipping, preserve liftability.
\end{lem}

\begin{proof}
This is well known and we omit its proof, see ~\cite[Section~4.2]{Z95}.
\end{proof}

The whole construction will consist  in successively applying these operations
in such a way that at each iteration the new cells are so small that 
their intersection pattern with lines can be controlled.

We call a set of planar points in \defn{general position} if 
no three are collinear.

\begin{lem}\label{lem:truncating}
 Let $S$ be a subset of the vertices of a subdivision of $\RR^2$ that are in general position.
 Then the vertices of $S$ can be truncated in such a way that no line intersects three
 of the newly created regions. 
\end{lem}
\begin{proof}%
 From the general position assumption there is some $\delta>0$ such that any line through 
 two points in $S$ stays at distance at least $\delta$ from any third point.
 Hence, there exists an $\varepsilon>0$ such that any line that goes through two points, each at distance at most
 $\varepsilon$ from a different point of~$S$, stays at distance at least $\varepsilon$ from the remaining points of $S$.
 The claim follows from the fact that the truncation regions can be arbitrarily small around the truncated points.
\end{proof}

\begin{lem}\label{lem:unzipping}
 Let $T$ be a triangular region of a Euclidean subdivision, $\ell$ a line 
 through one of its vertices $v$ that intersects the interior of $T$, and $\varepsilon>0$ a real. Then $T$ can be unzipped
 towards $v$ in such a way that for every line $\ell'$ that intersects at least three 
 regions of the spine, the angle between $\ell$ and $\ell'$ is at most $\varepsilon$.
 
 This can be done even when one forces the new vertices to be in general position
 with respect to a given point configuration.
\end{lem}
\begin{proof}%
 Start by stacking with a point $v'$ on $\ell$. Notice that the truncations can be made with very thin triangles,
 in such a way that the spine is sufficiently close to the edge $(v,v')$ in Hausdorff distance. If the
 pieces have a long enough diameter with respect to the distance of the spine to the edge, then any vector whose endpoints 
 belong two non-consecutive pieces of the spine will form a very small angle with $(v,v')$. In particular, 
 the line spanned by these points can be forced to be arbitrarily close to the line $\ell$.
 
 The last claim follows from the freedom in the choice of the truncation points 
 (the starting line $\ell$ might have to be
 perturbed before starting if the configuration has points on it).
\end{proof}

\subsection{The upper bound: completion of proof}\label{subsec:ub2}

We are ready for the proof of Theorem~\ref{thm:UB}. %
\begin{figure}[htpb]

\includegraphics[width=.55\linewidth]{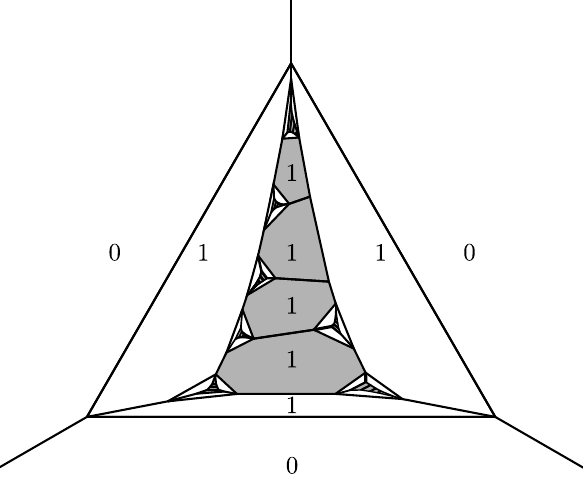}
 
 \caption{A schema of the construction in the proof of Theorem~\ref{thm:UB}, with $\ell=7$ and $k=2$. Numbers indicate the level of the regions (unnumbered regions are at level $2$), and spine regions are 
 shadowed.}\label{fig:proof_sketch}
\end{figure}

\begin{proof}[Proof of Theorem~\ref{thm:UB}]
A sketch of the  construction is depicted in Figure~\ref{fig:proof_sketch}.

The starting point of the construction is a regular simplex, inscribed on the unit sphere with one vertex 
at the south pole $(0,0,-1)$. We consider the polyhedral subdivision $\sE_0$ obtained by centrally projecting its spherical image subdivision.
It consists of a bounded triangular region $T_0$ and three unbounded regions. We say that these $4$
regions are at level~$0$. Note that $\sE_0$ is a liftable subdivision.

The triangle $T_0$ will be unzipped at length $t-3$, for some $t\ge5$ that will be defined later,
in such a way that all the points are in general position.
Then we will truncate $t$ of the $2t -5$ newly created vertices on the spine, in such a way that
that no line intersects three of the newly created regions, using Lemma~\ref{lem:truncating}. The new regions are at level~$1$ and $T_0$ is their
predecessor.

For~$i$ from $1$ to $k$ ($k$ will also be defined later), we will repeat this operation on all the triangles at level~$i$ (there are $t$ of them for each triangle at level $i-1$). This is done 
as follows. We process the triangles at level~$i$ one by one. First we select a 
line through one of its vertices
whose direction forms an angle of at least $2\varepsilon$ with all the lines chosen until now (in this and previous levels). 
This can be done
by choosing a set of well-separated candidate directions beforehand, one for each region that will have to be 
unzipped, and setting $\varepsilon$ accordingly.
We apply then Lemma~\ref{lem:unzipping} to unzip this triangle at length $t-3$ in such a way that any line through two of its non-consecutive spine regions
must form an angle of at most $\varepsilon$ with its line (and hence cannot intersect two non-consecutive spine regions of one of the previous spines); while keeping all new vertices in general position.

Except for the last iteration $i=k$, once this is done we choose $t$ among the new spine vertices in each triangle,
and we truncate them in such a way that no line intersects three of these newly created regions, using Lemma~\ref{lem:truncating}. These new triangular regions are at level~$i+1$ and their
predecessor is the triangle at level~$i$ that contained them.

Observe that, when unzipping, we replace each triangle at level~$i$ by $t$ new regions at level $i+1$ ($t-3$ of which are spine regions and $3$ non-spine regions). And then we create $t$ 
triangles at level $i+1$ by truncating the spine vertices (when $i<k$).
This way, the number of regions at level~$i$ is $3$ for $i=0$ and $t^i$ for $1\leq i\leq k$. That is, the total number of regions is 
$$n=2+\frac{t^{k+1}-1}{t-1},$$
and therefore $k\leq \log_t(n)$.

We compute now the maximal number of regions that can be intersected by a line. 
By construction, if a line intersects more than $2$ regions of a spine, then it cannot intersect more than two regions from any other spine. Hence, except for maybe one
spine where it can go through at most $t-3=\sO(t)$ regions, it intersects at most $2$ regions from the remaining spines. We count these $\sO(t)$ separately and
continue counting as if no line could intersect more than $2$ regions of any spine.

Hence, for a triangle at level~$i$, a line can intersect at most $3$ non-spine regions and $2$ spine regions at level~$i+1$. Thus, for each triangle, 
there are at most $5$ regions that have it as predecessor that intersect any given line. 
For each level $i\geq 1$, no line can intersect more than two triangles at level~$i$ (because we used Lemma~\ref{lem:truncating}). Since there are $k$ levels $\geq 1$, 
this amounts for at most $10\cdot k$ regions intersected by any single line. And there are at most $3$ regions at level $0$.
These are $\sO(k)$ regions that can be intersected in addition to the at most $\sO(t)$ regions in a single spine. Hence, a line crosses at most
$\sO(t+k)=\sO(t+\log_{t}(n))$ regions.

Taking $t=\left\lfloor\frac{\log(n)}{\log(\log(n))}\right\rfloor$ gives that at most
$$\sO\left(\frac{\log(n)}{\log(\log(n))}\right)$$
regions are intersected by any line. Note that any large enough value of~$n$ can be attained by this construction just by taking $t=\left\lfloor{\log(n)}/{\log(\log(n))}\right\rfloor$, 
$k=\lceil\log_t (n)\rceil$, and adjusting the length at which the triangles are unzipped at the last iteration.

Since all the operations were liftable by Lemma~\ref{lem:liftable}, we can lift this subdivision to the spherical subdivision corresponding to a polytope $P$ with $n$ vertices. 
Since there are only three regions of $\sDD(P)$ intersecting the upper hemisphere, the great circle span of $\sDD(P)$ and that of its intersection with the lower hemisphere differ at most by three (see
Lemma~\ref{lem:equivcirclelinespan}).
Therefore, by Lemma~\ref{lem:equivalencegreatcircleshadow}, the shadow number of $P$ is at most
\[\sh(P)=\sO\left(\frac{\log(n)}{\log(\log(n))}\right).\qedhere\]
\end{proof}

\section{Moser's unbounded shadow problem}\label{sec:unbounded}

In this section, we will determine the shadow number for unbounded polyhedra (whose size is measured in terms of their number $n$ of bounded vertices).
There are two results. In Proposition \ref{prop:lower} we give a lower bound showing $\sh(P_n) \ge 3$
for $n \ge 3$.
In Theorem~\ref{thm:polyhedrality} we will construct a sequence of unbounded polyhedra $P_n$, for all $n\geq 4$, having $n$ vertices 
and $n$ faces and whose shadow number is $\sh(P_n) = 3$, giving an upper bound for $\shf_u(n)$.\footnote{We are grateful to an anonymous reviewer who suggested this example to improve our original upper bound of~$5$ to~$3$.}
 Both results together establish Theorem~\ref{thm:UnboundedMoserShadowProblem}.

\subsection{Unbounded shadow lower bound}\label{subsec:unbounded-shadow-lower}

The following proposition gives a lower bound for the unbounded shadow number function. 
\begin{prop}\label{prop:lower}
Every unbounded polyhedron $P$ with at least $3$ vertices has shadow number greater or equal to~$3$ (and at least $n$ for $n\leq 3$). In particular, $\shf_u(n) \geq 3$ for all $n\geq 3$ (and $\shf_u(1) = 1$ and $\shf_u(2) = 2$). \end{prop}
\begin{proof}
Given $P$, let  $\sD=\sDD(P)$ be the spherical image subdivision that $P$ induces on $U = \sigma(P)$. Since $P$ is unbounded, $U$ lies in a closed hemisphere. Without loss of generality we may  assume that the interior of $U$ lies in $\SSS_{-}^2$. 
The central projection of~$\sD$ is a Euclidean subdivision~$\sE$ of a polygon~$Q$.
 Since~$P$ has at least three vertices, $\sE$ consists of at least three regions, $U_1$, $U_2$ and $U_3$. 
 We can assume without loss of generality that $U_1$ and $U_2$ share an edge. 
 Take $x$ to be a point in the relative interior of this edge, and $y$ a point in the interior of~$U_3$. 
 After a small perturbation if needed, the line through $x$ and $y$ stabs the interior of $U_1$, $U_2$ and $U_3$. This fact implies that $\shf(P) \ge 3$  by Lemma~\ref{lem:equivcirclelinespan}. 
For two or  fewer  regions, there is always a line in $\sE$ through all of them.
\end{proof}

\subsection{Unbounded shadow problem: upper bound}\label{subsec:unbounded-shad--upper}

The following construction gives an upper bound for unbounded shadow number function. 

\begin{theorem}\label{thm:polyhedrality}
For each $n \ge 4$, there is an unbounded pointed convex polyhedron with $n$ vertices~$P_n$ whose shadow number $\sha(P_n)$ is~$3$; consequently $\shf_u(n) \leq 3$ for all $n$ (it is trivially true for $n\leq 3$). %
\end{theorem}

\begin{proof}
 For $n\ge 4$, consider the convex polyhedral cone 
 \[Q_n:=\left\{x\in\RR^3\ : \: \sprod{x}{w_k}\leq 0, \text{ for }  0\leq k\leq n-1\right\},\]
where $w_k:=\left( \cos\left(\frac{2\pi k}{n-1}\right) , \sin\left(\frac{2\pi k}{n-1}\right), -1\right)$ and $\sprod{\cdot}{\cdot}$ denotes the standard scalar product. This is a cone over a regular $(n-1)$-gon. It has a single vertex at the origin and $n-1$ (unbounded) facets. Now we stack a vertex on top of each of these $n-1$ facets. That is, for each facet we add a point that is slightly beyond it and beneath the hyperplanes defining the remaining facets, and take the convex hull. We obtain an unbounded polyhedron with~$n$ vertices (the origin plus $n-1$ stacking points), and $3(n-1)$ (unbounded) facets (see Figure~\ref{fig:UnboundedShadow}, left). One explicit realization is the following polyhedron~$P_n$:
 \begin{align*}
  P_n:=
  \big\{
  x\in\RR^3\ : \ & \sprod{x}{w_k}\leq 1,\  \sprod{x}{\tfrac{2}{3}w_k+\tfrac{1}{3}w_{k+1}}\leq 0, \text{ and }\\
			    & \sprod{x}{\tfrac{2}{3}w_k+\tfrac{1}{3}w_{k-1}}\leq 0; \text{ for }  0\leq k\leq n-1
  				    \big\}.
 \end{align*}

\begin{figure}[htpb]
\centering
\includegraphics[width=.25\linewidth]{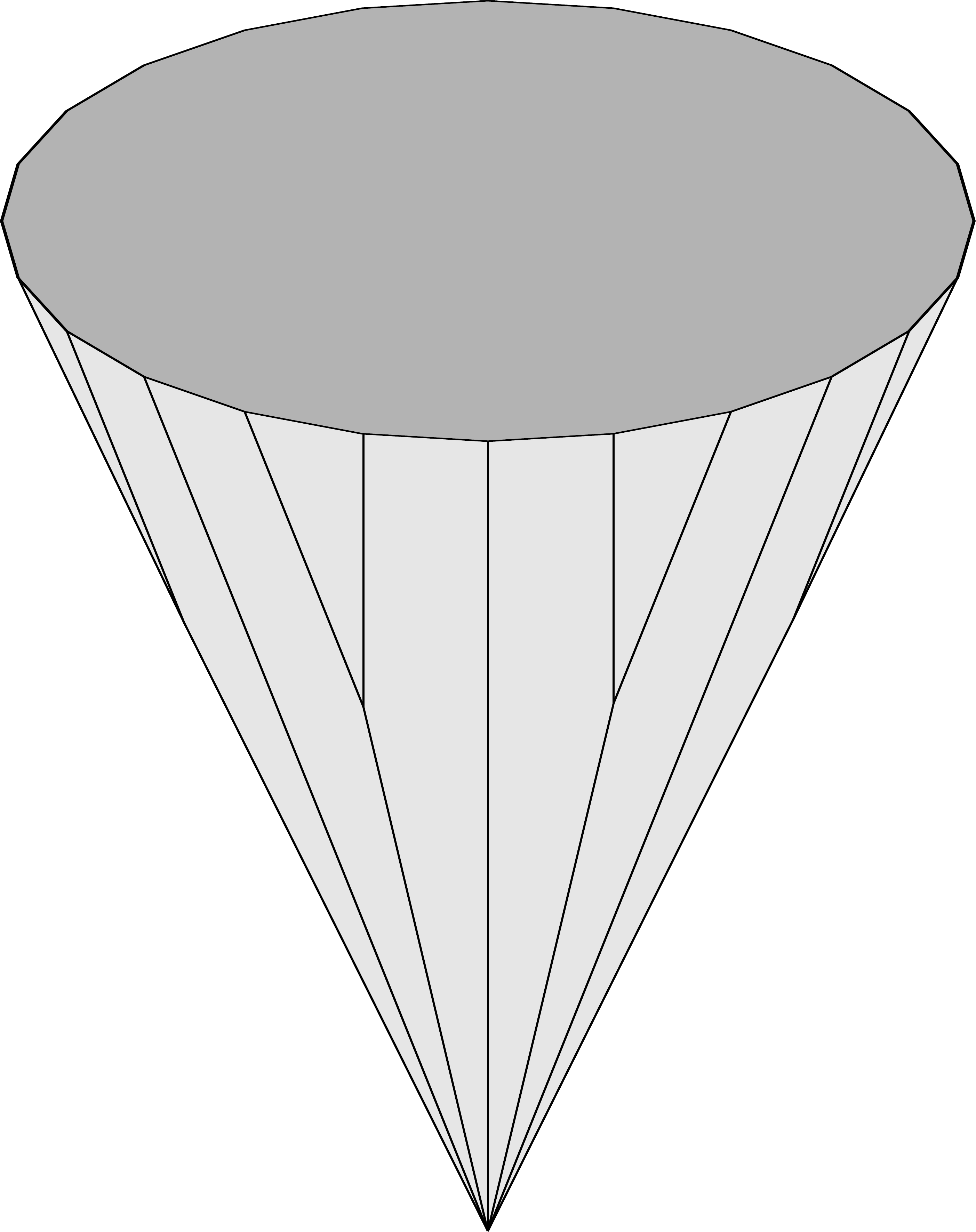}
\qquad\qquad\qquad
\raisebox{.25cm}{\includegraphics[width=.27\linewidth]{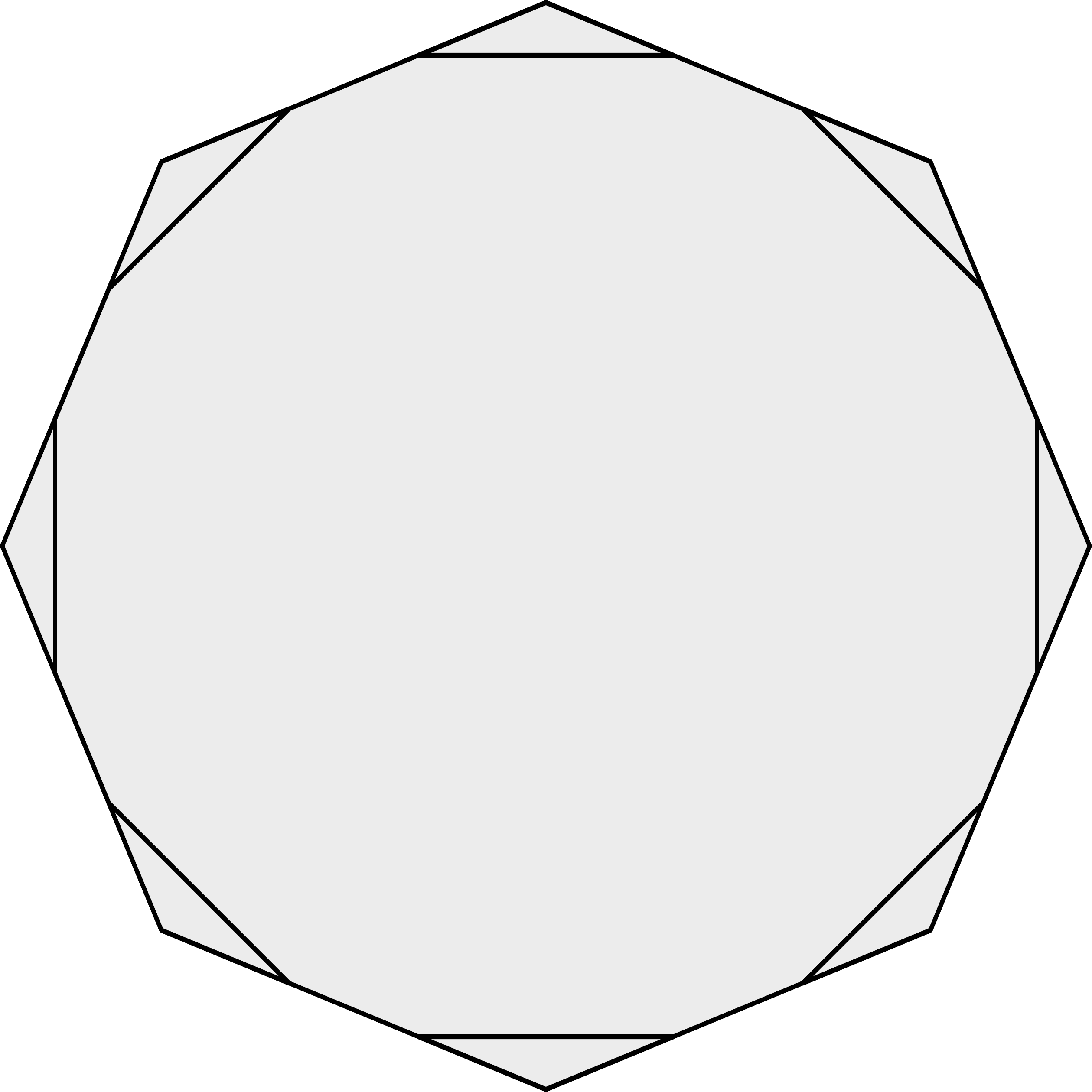}}
\caption{An instance of $P_n$, for $n=9$, and the central projection of its spherical image subdivision, which is an $8$-gon subdivided into 9 regions.}
\label{fig:UnboundedShadow}
\end{figure}
The shadow number of $P_n$ is at most~$3$. Indeed, let $\pi:\RR^3\to \RR^2$ be a linear projection. If~$\pi(Q_n)$ does not cover the whole plane then it is a two-dimensional cone pointed at the origin and bounded by the image of two of the rays of~$Q_n$, which are also rays of~$P_n$. Besides the origin, only the vertices of $P_n$ stacked to facets incident to these rays can appear as vertices of the shadow~$\pi(P_n)$. Moreover, for each of the two sides, only one of the two neighboring stacked vertices can be visible: They cannot be both separated from the shadow of~$Q_n$, as the segment between both intersects the interior of~$Q_n$ (because each ray of~$Q_n$ is preserved by the stacking operation).

This can also be easily seen on its spherical image subdivision~$\sD_{\sigma}(P_n)$. The spherical image subdivision of~$Q_n$ is a spherical regular $(n-1)$-gon centered at the south pole. The stacking operation will subdivide each of its edges into three pieces and add triangles joining each original vertex with its two neighboring newly created vertices (as in Figure~\ref{fig:UnboundedShadow}, right). It is straightforward to check that the great circle span of this subdivision is at most~$3$.
\end{proof}

\begin{rmk}\label{rmk:size2}
In connection with Remark~\ref{rmk:size},
the  examples $P_n$ used  in the construction of Theorem~\ref{thm:polyhedrality} work for the edge and facet numbers size measures as well.
These $P_n$ have $n$ vertices, $3n-3$ facets, and  $4n-4$ edges, and so establish boundedness of
 the shadow function with respect to  any of these measures of size of the polyhedron going to $+\infty$.
\end{rmk}

\begin{rmk}
In contrast with the shadow number, the silhouette span of these polyhedra is $\sh^{\ast}(P_n) \ge 2(n-1)$. This can be seen
by taking a point $p$ just above the apex of the pointed cone $P_n$, since the $n-1$ stacked vertices and the $n-1$ rays of $Q_n$ are visible as vertices of the silhouette of $P_n$ viewed from~$p$.
The stacked vertices would be visible from this point but the $n-1$ rays would be missing in the pre-silhouette definition 
in \cite[Section 5.3, p.174]{CEG89}.
\end{rmk}

\section{The unbounded silhouette span problem}\label{sec:unboundedsilhouette}

In this final section, we consider silhouettes of possibly unbounded polyhedra, and 
determine the asymptotics of  the unbounded silhouette span function $\siu(n)$.

We will show that the asymptotic growth rate of $\sih(n)$ and $\sib(n)$ are of the same
order by reducing the unbounded case to the bounded case using a projective transformation.

\begin{theorem}
The unbounded $n$-vertex silhouette span number~$\siu(n)$ for $3$-dimensional convex polyhedra satisfies  
\begin{equation*}
                                                   \siu(n) = \Theta\left(\frac{\log (n)}{\log(\log (n))}\right).
\end{equation*}
\end{theorem}

\begin{proof}
The upper bound follows from the trivial inequality $\siu(n)\leq \sib(n)$ and Theorem~\ref{thm:BoundedSilhouette}.
It suffices thus to prove the lower bound
\begin{equation*}
                                                   \siu(n) = \Omega\left(\frac{\log (n)}{\log(\log (n))}\right).
\end{equation*}
This lower bound holds for polytopes by Theorem~\ref{thm:BoundedSilhouette}, so we concentrate on unbounded polyhedra. 

Let~$P$ be an unbounded polyhedron with $n>0$ vertices (which is therefore pointed). The spherical image subdivision of $P$ is a spherical polygon contained in a closed hemisphere. After a suitable rotation, we may assume that this spherical polygon lies in the lower hemisphere and contains the south pole. This implies that every vertex admits a supporting hyperplane with a normal vector
that has a negative third coordinate, and that for $M\in\R$ large enough the hyperplane $H_M$, defined by $H_M=\{(x,y,z)\ :\ z=-M\}$, avoids $P$.

We will take some very large $M\gg 0$ with $H_M\cap P=\emptyset$ and consider the following projective transformation $\phi$ sending $H_M$ to infinity:
\begin{align*}
\phi\ :\ \R^3 \smallsetminus H_M & \to  \R^3 \smallsetminus H_{-M}\\ (x,y,z) & \mapsto \frac{(x,y,z)}{1+\frac{z}{M}}.
\end{align*}
It maps bijectively $\R^3 \smallsetminus H_M$ to $\R^3 \smallsetminus H_{-M}$, with inverse $\phi^{-1}(x,y,z)={(x,y,z)}/{(1-\frac{z}{M})}$, and sends affine subspaces to affine subspaces, preserving incidences (see~\cite[Appendix~2.6]{Z95} for a brief introduction to projective transformations in the context of polyhedra). In particular, it maps hyperplanes to hyperplanes as follows:
\begin{equation}
\{(x,y,z)\ :\ ax+by+cz+d=0\} \mapsto \{(x,y,z)\ :\ ax+by+(c-\frac{d}{M})z+d=0\}.\label{eq:hyperplanemap}
\end{equation}
The closure of the image of the polyhedron $P$ is the (bounded) polytope $Q$ bounded by the inequalities inherited from $P$ via~\eqref{eq:hyperplanemap} together with the new inequality $z\leq M$. Let $F$ be the face of $Q$ supported by the hyperplane $\{z=M\}$. Then $\phi(P)=Q\setminus F$, as $F$ is the ``face at infinity'' of~$P$.

Now, since every vertex of~$P$ has a supporting hyperplane pointing downwards, so do the vertices of~$Q$ that do not belong to $F$ (provided that $M$ is large enough). Therefore, the restriction of the spherical image subdivision of $Q$ to the lower hemisphere has at least $n$~regions. 
We can now follow the proof of the lower bound Theorem~\ref{thm:LB} and find a direction $\bv=(v_1,v_2,v_3)$ such that the shadow of~$Q$ in direction $\bv$ has at least 
$\Omega\left(\frac{\log (n)}{\log(\log (n))}\right)$ vertices of $Q\smallsetminus F$. 
Since small perturbations do not decrease the shadow number, we can assume that $v_3\neq 0$.

When $v_3\neq 0$, lines in direction $\bv$ are mapped by $\phi^{-1}$ to lines through the point $p=(\frac{-Mv_1}{v_3}, \frac{-Mv_2}{v_3}, -M)$. Consequently, lines in direction $\bv$ through a point 
of $\phi(P)=Q\setminus F$ are mapped by $\phi^{-1}$ bijectively to lines through~$p$ and a point in~$P$. 
In fact, $\phi(C_p(P)\smallsetminus p)$ is easily seen to be the one-sided cylinder 
$(Q+\bv\RR) \cap \{z\leq M\}$. Hence, the shadow of $Q$ in direction~$\bv$ has the same number of vertices
 as the silhouettes of $P$ from $p$. This can be seen explicitly by noting than the image of 
 the silhouette $C_p(P)\cap H_{\frac{M}{2}}$ under the projective transformation $\phi$ is
  the polygon $(Q+\bv\RR)\cap H_{{M}}$, together with the fact that (admissible) projective transformations
   do not change the combinatorial type.

Thus, the silhouette span of $P$ is $\Omega\left(\frac{\log (n)}{\log(\log (n))}\right)$.
\end{proof}

\subsection*{Acknowledgments}
We are indebted to anonymous reviewers whose remarks helped improve the paper; and in particular for suggesting the example in Section~\ref{sec:unbounded}, which allowed for improving our original upper bound for the unbounded shadow number from~$5$ to~$3$.
Some work of the second author was
done while at the National University of Singapore.
He thanks them for support to work with the first author
in Michigan.

\end{document}